\documentclass{amsart}

\begin{document}
\title[Hanna Neumann Conjecture]{A short proof that positive
  generation implies the Hanna Neumann Conjecture} \author{Walter D.
  Neumann} \address{Department of Mathematics\\Barnard College,
  Columbia University\\New York, NY 10027, USA}
\email{neumann@math.columbia.edu}
\maketitle

The Strengthened Hanna Neumann Conjecture \cite{neumann} posits that if
$U$ and $V$ are finitely generated subgroups of a free group $F$ then
$$\sum_{x\in T}\chi_0(U\cap x^{-1}Vx)\le \chi_0(U)\chi_0(V)\,,$$
where $\chi_0(H)=\operatorname{max}(0, \operatorname{rank}\, H - 1)$
and $T$ is a set of double coset representatives for $V\backslash
F/U$. This is proved in \cite{khan} and \cite{meakin_weil} when $U$ is
positively generated (i.e., generated by elements of the subsemigroup
generated by a basis of $F$). In a February 2003 CUNY Group Theory
Seminar I pointed out a simple proof. It was suggested then and
again recently that I record this in print.

Since any free group embeds in the free group of rank 2, we may assume
$F=\langle a,b\rangle$. By embedding $F$ in itself using the map
$a\mapsto a^2$, $b\mapsto ab$, we may further assume that $U$ is
generated by positive words in the elements $a^2$ and $ab$. Then the
core graph $G_0(U)$ has just two types of valence 3 vertices (since
this is true of $G_0(\langle a^2, ab\rangle)$) and they occur in
exactly equal numbers (since they have (incoming, outgoing) valence
$(2,1)$ and $(1,2)$ respectively and the strings of $G_0(U)$ are
directed).  Proposition 3.1 of \cite{neumann} says that a
counterexample to the generalized HN conjecture must have over half
the valence 3 vertices of both $G_0(U)$ and $G_0(V)$ all of the same
type, so $U$ can't belong to a counterexample.


\begin{thebibliography}{BDM}

\bibitem{khan} Bilal Khan.
Positively generated subgroups of free groups and the Hanna Neumann
conjecture. {\it Combinatorial and geometric group
theory (New York, 2000/Hoboken, NJ, 2001)}, 155--170,
Contemp. Math., 296, Amer. Math. Soc., Providence, RI, 2002.


\bibitem{meakin_weil} J. Meakin, P. Weil.
Subgroups of free groups: a contribution to the Hanna Neumann
conjecture.
Proceedings of the Conference on Geometric and Combinatorial Group
Theory, Part I (Haifa, 2000).
Geom. Dedicata 94 (2002), 33--43.

\bibitem{neumann} Walter D. Neumann. On intersection of finitely
  generated subgroups of free groups. {\it Groups -- Canberra 1989.}
  Lecture Notes in Mathematics {\bf 1456}, 161--170 (Springer Verlag,
  Berlin, Heidelberg, New York, 1990). 
\end{thebibliography}
\end{document}